\ifx\documentclass\undefined
\documentstyle[12pt]{article}
\else
\documentclass[12pt]{article}
\usepackage{latexsym}
\usepackage{amsfonts}
\usepackage{amsmath}
\usepackage{amssymb}
\fi

\author{K. Bezdek \thanks{Partially supported by a Natural Sciences and 
Engineering Research Council of Canada Discovery Grant.} (Univ. of Calgary) and
\and Gy. Kiss  \thanks{Partially supported by
the South African-Hungarian Intergovernmental
Scientific and Technological Cooperation Project, Grant No.
ZA-21/2006, by 
the Hungarian National
Foundation for Scientific Research, Grant No. NK 67867 and by
the Hungarian National Office of Research and Technology within the framework of 
``\"Oveges J\'ozsef'' program. } (E\"otv\"os Univ.)}

\font\tenBbb=msbm10 at 12pt         \font\sevenBbb=msbm9    \font\fiveBbb=msbm7
\newfam\Bbbfam
\textfont\Bbbfam=\tenBbb \scriptfont\Bbbfam=\sevenBbb
\scriptscriptfont\Bbbfam=\fiveBbb

\def\kkk{\null\hfill $\Box$\smallskip}

\def\KK{{\bf K}}

\def\PP{{\bf P}}

\def\CC{{\bf C}}
\def\MM{{\bf M}}

\def\WW{{\bf W}}

\newcommand{\proof}{{\noindent\bf Proof:{\ \ }}}

\newtheorem{theorem}{Theorem}[section]
\newtheorem{lemma}[theorem]{Lemma}

\newtheorem{proposition}[theorem]{Proposition}

\newtheorem{cor}[theorem]{Corollary}
\newtheorem{con}[theorem]{Conjecture}

\newtheorem{definition}[theorem]{Definition}

\title{On the X-ray number of almost smooth convex bodies and of convex bodies of constant width
}

\begin{document}

\maketitle

\date

\begin{abstract}
The X-ray numbers of some classes of convex bodies are investigated.
In particular, we give a proof of the X-ray Conjecture as well as of the Illumination Conjecture for almost smooth convex bodies of any dimension and for convex bodies of constant width of dimensions $3$, $4$, $5$ and $6$.
\end{abstract}

\section{Introduction}
\label{zero}
In 1972, the X-ray number of convex bodies was introduced by P. Soltan as follows (see also \cite{ms}). Let $\KK$ be a convex body of ${\bf E}^d, d\ge 2$ (i.e. a
compact convex set of the $d$-dimensional Euclidean space ${\bf
E}^d$ with non-empty interior). Let $L\subset {\bf E}^d$ be a line
through the origin. We say that the point
$p\in {\KK}$ is {\it X-rayed along $L$} if the line parallel to $L$
passing through $p$ intersects the interior of ${\KK}.$ The
{\it X-ray number} $X({\KK})$ of ${\KK}$ is the smallest number of
lines such that every point of ${\KK}$ is X-rayed along at least one
of these lines. Obviously, $X({\KK})\geq d$. Moreover, it is easy to see that
this bound is attained by any smooth convex body. On the other hand, if
${\CC}_d$ is a $d$-dimensional (affine) cube and $F$ is one of its $(d-2)$-dimensional 
faces, then the X-ray number of the convex hull of the set of vertices
of ${\CC}_d\setminus F$ is $3\cdot 2^{d-2}.$

In 1994, the first named author of this note and Zamfirescu \cite{bz} published the conjecture that the X-ray
number of any convex body in ${\bf E}^d$ is at most $3\cdot 2^{d-2}$.
This conjecture, which we call the X-ray Conjecture, is proved only in the plane and it is open in high dimensions. A related and much better studied
problem is the Illumination Conjecture of Boltyanski and Hadwiger according to which any $d$-dimensional convex body can be illuminated by $2^d$ directions (resp., point sources). For a recent account on the status of the Illumination Conjecture we refer the interested reader to \cite{bk06} as well as to \cite{ms}.
Here we note that if $I(\KK)$ denotes the minimum number of directions that are needed for the illumination of the convex body $\KK\subset {\bf E}^d$, then the inequalities $X(\KK)\le I(\KK)\le 2\cdot X(\KK)$ hold.
Putting it differently, any proper progress on the X-ray Conjecture would imply a progress on Illumination Conjecture and vica versa. Last but not least we note that a natural way to prove the X-ray Conjecture would be to show that any convex body ${\KK}\subset {\bf E}^d$ can be illuminated by $3\cdot 2^{d-2}$ pairs of pairwise opposite directions.

The main goal of this paper is to study the X-ray numbers of almost smooth convex bodies and of convex bodies of constant width. As a result we get a proof of the X-ray Conjecture as well as of the Illumination Conjecture for almost smooth convex bodies of any dimension and for convex bodies of constant width in dimensions $3$, $4$, $5$ and $6$. It would be interesting to extend the method of this paper for the next couple of dimensions (more exactly, for dimensions $7\le d\le 15$) in particular, because in 
these dimensions neither the X-ray Conjecture nor the Illumination Conjecture are known to hold for convex bodies of constant width (for more details see \cite{schr}). Another aspect of this paper, is to encourage further research on covering a $(d-1)$-dimensional unit sphere of ${\bf E}^d$ with a given number of pairwise antipodal congruent spherical caps of smallest possible spherical radius. Namely, based on Lemma~\ref{alt} any properly chosen improvement on the above covering problem (in particular, for the values of $d$ mentioned above) may lead to an improved upper bound on the X-ray number of convex bodies of constant width in ${\bf E}^d$. Also, as it is discussed in the last section of this note, the X-ray Conjecture implicitly contains a conjecture that is strongly connected to the elegant theorem of Danzer and Gr\"unbaum \cite{dg} on antipodal convex polytopes.

\section{The X-ray number of almost smooth convex bodies}
\label{one}

\begin{definition}
Let ${\MM}\subset {\bf E}^d$ be a convex body and $b$ be a boundary point
of $\MM$ (i.e. $b\in {\rm bd}\  {\MM}$). Let $N_b$ denote the set of outer normal vectors of
all supporting hyperplanes of $\MM$ at $b$. Then $\MM$ is called almost
smooth, if for each $b\in {\rm bd}\  {\bf M}$ and any ${\bf n}_i, {\bf n}_j\in N_b$
the inequality
$$\frac{ {\bf n}_i\cdot {\bf n}_j}{ \| {\bf n}_i\| \cdot \| {\bf n}_j \| }
\geq \frac{d-2}{d-1}$$
holds, where $\ \cdot\ $ in the numerator refers to the standard inner product of ${\bf E}^d$.
\end{definition}

Let $\KK$ be a convex body in ${\bf E}^d$ and let $F$ be a face of
${\KK}$ that is let $F$ be the intersection of $\KK$ with some of its supporting hyperplanes. 
The {\it Gauss image} $\nu (F)$ of the face $F$ is the set of
all points (i.e. vectors) $n$ of the $(d-1)$-dimensional unit sphere $S^{d-1}\subset
{\bf E}^d$ centered at the origin $o$ of ${\bf E}^d$ for which the supporting
hyperplane of ${\KK}$ with outer normal vector $n$ contains $F.$
It is easy to see that the Gauss images of distinct faces of $\KK$
have disjoint relative interiors and $\nu (F)$ is compact and spherically
convex for any face $F.$

 Recall the following statement published in \cite{bz} that gives a reformulation of the X-ray number of a convex body in terms of its Gauss map.

\begin{lemma} \label{elso}
Let $\KK$ be a convex body in ${\bf E}^d, \, d>2,$ and let $b\in {\rm bd}\ {\KK}$
be given moreover, let $F$ denote (any of) the face(s) of $\KK$ of smallest dimension containing $b$.
Then $b$ is X-rayed along the line $L$ if and only if
$L^{\perp }\cap \nu(F)=\emptyset$, where  $L^{\perp }$ denotes the hyperplane orthogonal to $L$ and passing through the origin $o$ of ${\bf E}^d$. Moreover, $X({\KK})$ is the smallest
number of $(d-2)$-dimensional great spheres of $S^{d-1}$ with the property that the
Gauss image of each face of $\KK$ is disjoint from at least one of
the given great spheres.
\end{lemma}

As an easy application one can obtain the following statement. If $\MM$ is a smooth convex body in ${\bf E}^d, \, d>1$ (i.e. each boundary point of $\MM$ belongs to exactly one supporting hyperplane), then $X({\MM})=d$. In fact, the following stronger result holds.

\begin{theorem} \label{masodik}
If ${\MM}\subset {\bf E}^d,$ $d\geq 3$ is an almost smooth convex body, then
$$X({\MM})=d.$$
\end{theorem}

\proof Recall the following notion. Let $C\subset S^{d-1}$ be a set of finitely many points. Then the {\it covering radius} of $C$ is the smallest positive real number $r$ with the property that the family of spherical balls of radii $r$ centered at the points of $C$ cover $S^{d-1}$. Our proof of Theorem~\ref{masodik} relies on the following rather general concept.

\begin{lemma} \label{harmadik}
Let ${\MM}\subset {\bf E}^d$, $d\geq 3$ be a convex body and let $r$ be a positive real number with the property that the Gauss image $\nu (F)$ of any face $F$ of ${\MM}$ can be covered by a spherical ball of radius $r$ in $S^{d-1}$. Moreover, assume that there exist $2m$ pairwise antipodal points of $S^{d-1}$ with covering radius $R$ satisfying the inequality $r+R\le\frac{\pi}{2}$. Then $X(\MM)\le m$. 
\end{lemma}

\proof
Let $\{p_1,-p_1, \dots, p_m,-p_m\}$ be the family of pairwise antipodal points in $S^{d-1}$ with covering radius $R$. Moreover, let $B_i\subset S^{d-1}$  be the union of the two $(d-1)$-dimensional closed spherical balls of radius $R$ centered at the points $p_i$ and $-p_i$ in $S^{d-1}$, $1\le i\le m$. Finally, let $S_i$ be the $(d-2)$-dimensional great sphere of $S^{d-1}$ whose hyperplane is orthogonal to the diameter of $S^{d-1}$ with endpoints $p_i$ and $-p_i$, $1\le i\le m$. Based on Lemma~\ref{elso} it is sufficient to show that the
Gauss image of each face of $\MM$ is disjoint from at least one of
the great spheres $S_i, 1\le i\le m$. 

Now, let $F$ be an arbitrary face of the convex body ${\MM}\subset {\bf E}^d$, $d\geq 3$, and let $B_F$ denote the smallest spherical ball of $S^{d-1}$ with center $f\in S^{d-1}$ which contains the Gauss image $\nu (F)$ of $F$. By assumption the radius of $B_F$ is at most $r$. As the family $\{B_i, 1\le i\le m\}$ of antipodal pairs of balls forms a covering of $S^{d-1}$ therefore $f\in B_j$ for some $1\le j\le m$. If in addition, we have that $f\in {\rm int} B_j$ (where ${\rm int}(\ )$ denotes the interior of the corresponding set in $S^{d-1}$), then the inequality $r+R\le\frac{\pi}{2}$ implies that $\nu (F)\cap S_j=\emptyset$. If $f$ does not belong to the interior of any of the sets $B_i, 1\le i\le m$, then clearly $f$ must be a boundary point of at least $d$ sets of the family $\{B_i, 1\le i\le m\}$. Then either we find an $S_i$ being disjoint from $\nu (F)$ or we end up with $d$ members of the family $\{S_i, 1\le i\le m\}$ each being tangent to $B_F$ at some point of $\nu (F)$. Clearly, the later case can occur only for finitely many $\nu (F)$'s and so, by taking a proper congruent copy of the great spheres $\{S_i, 1\le i\le m\}$ within $S^{d-1}$ (under which we mean to avoid finitely many so-called prohibited positions) we get that each 
$\nu (F)$ is disjoint from at least one member of the family $\{S_i, 1\le i\le m\}$. This completes the proof of Lemma~\ref{harmadik}.
\kkk

First, note that according to the spherical version of the well-known Jung theorem \cite{d} the Gauss image $\nu (F)$ of any face $F$ of the almost smooth convex body ${\MM}\subset {\bf E}^d,$ $d\geq 3$ can be covered by a spherical ball of radius $r=\arccos \sqrt{\frac{d-1}{d}}$ in $S^{d-1}$ where, $r$ is in fact, equal to the circumradius of a regular $(d-1)$-dimensional spherical simplex of edge length $\arccos(\frac{d-2}{d-1})$.

Second, take the $2d$ vertices of an arbitrary cross-polytope inscribed in $S^{d-1}$. An easy computation shows that the covering radius of the $2d$ pairwise antipodal points of $S^{d-1}$ just introduced is $R=\arccos \sqrt{\frac{1}{d}}$.

Finally, as $r+R=\frac{\pi}{2}$ therefore Lemma~\ref{harmadik} finishes the proof of Theorem~\ref{masodik} in a straightforward way.
\kkk

\section{On the X-ray number of convex bodies of constant width}
\label{two}

Let $\WW\subset {\bf E}^d$ be a convex body of constant width $1$.
Moreover, let $p$ and $q$ be arbitrary points of $\bf W$. Then $\| p-q\| \leq 1$
and $\| p-q\| =1$ if and only if there are parallel supporting hyperplanes
of $\WW$ at $p$  and $q$. Hence, if $p\in {\rm bd}{\WW}$ and
$N_p$ is the set of outer normal vectors of
the supporting hyperplanes of $\WW$ at $p$, then for all
${\bf n}_i, {\bf n}_j\in N_p$
the inequality
$$\frac{ {\bf n}_i\cdot {\bf n}_j}{ \| {\bf n}_i\| \cdot \| {\bf n}_j \| }
\geq \cos \frac{\pi }{3}$$
holds, otherwise the distance between the endpoints of the
inner unit normals
${\bf n}_i/\| {\bf n}_i \| $ and ${\bf n}_j /\| {\bf n}_j\| $
would be greater than $1$, a contradiction (because the diameter of $\WW$ is equal to $1$). Schramm \cite{schr} managed to combine this simple fact with some clever random techniques 
and proved the inequality $I(\WW)< 5d\sqrt{d}(4+\ln d)\big(\frac{3}{2}\big)^{\frac{d}{2}}$. Thus, also the inequality $X(\WW)< 5d\sqrt{d}(4+\ln d)\big(\frac{3}{2}\big)^{\frac{d}{2}}$ holds for any convex body of constant width ${\WW}\subset {\bf E}^d$ and for all $d\geq 3$. This inequality can be improved in small dimensions as follows.

First, note that any convex body of constant width $\WW$ in ${\bf E}^3$ is in fact, an almost smooth convex body of ${\bf E}^3$ and therefore Theorem~\ref{masodik} implies that in this case $X(\WW)=3$ and so, $I(\WW)\le6$. The later inequality is not new namely, it has been independently proved by quite a number of people (see the paper \cite{be07} for more details).   

Second, note that the following statement is a straightforward corollary of Lemma~\ref{harmadik}.

\begin{lemma}
\label{alt}
If there are $2m$
pairwise antipodal points on $S^{d-1}$ with covering radius
$r$ satisfying the inequality $r\leq \pi /2-r_{d-1}$,where $r_{d-1}=\arccos \sqrt{\frac{d+1}{2d}}$ is the circumradius of a 
regular $(d-1)$-dimensional spherical simplex of edge length $\pi /3$, then $X(\WW)\le m$ holds for any
convex body of constant width $\WW$ in ${\bf E}^d$. 
\end{lemma}

Now, we are ready to prove the X-ray Conjecture for convex bodies of constant width of dimensions $4,5$ and $6$.

\begin{theorem} \label{negyedik}
If $\WW$ is a convex body of constant width in ${\bf E}^4$, then $X(\WW)\leq 6$. Moreover, if $\WW$ is a convex body of constant width in ${\bf E}^d$ with $d=5, 6$, then $X({\WW})\leq 2^{d-1}$. 
\end{theorem}

\proof
Based on Lemma~\ref{alt} it is sufficient to find $12$ pairwise antipodal points of $S^3$ whose covering radius
is at most $\alpha=\pi /2-r_3=\pi /2-\arccos \sqrt{\frac{5}{8}}$. In order to achieve this let us take two regular hexagons of edge length $1$ inscribed into $S^3$ such that their $2$-dimensional planes are totally orthogonal to each other in ${\bf E}^4$. Now, let $\PP$ be the convex hull of the $12$ vertices of the two regular hexagons. If $F$ is any facet of $\PP$, then it is easy to see that $F$ is a $3$-dimensional simplex having two pairs of vertices belonging to different hexagons with the property that each pair is in fact, a pair of two consecutive vertices of the relevant hexagon. As an obvious corollary of this we get that if one projects any facet of $\PP$ from the center $o$ of $S^3$ onto $S^3$, then the projection is a $3$-dimensional spherical simplex whose two opposite edges are of length $\frac{\pi}{3}$ and the other four remaining edges are of length $\frac{\pi}{2}$. Also, it is easy to show that the circumradius of that spherical simplex is equal to $\alpha=\arccos \sqrt{\frac{3}{8}}$. This means that the covering radius of the $12$ points in question lying in $S^3$ is precisely $\alpha$ finishing the proof of Theorem~\ref{negyedik} for $d=4$.

In dimensions $d=5, 6$ we proceed similarly using Lemma~\ref{alt}. More exactly, we are going 
to construct $2^d$ pairwise antipodal points on $S^{d-1}$ ($d=5, 6$) with
covering radius at most 
$r\leq \pi /2-\arccos \sqrt{\frac{d+1}{2d}} =\arccos \sqrt{\frac{d-1}{2d}}.$

For $d=5$ we need to find $32$ pairwise antipodal points on $S^4$ with
covering radius at most $\arccos \sqrt{\frac{2}{5}}=50.768...^{\circ }$. 
Let us take a $2$-dimensional plane ${\bf E}^2$ and a $3$-dimensional subspace
${\bf E}^3$ in ${\bf E}^5$ such that they are totally orthogonal to each
other. Let ${\PP}_2$ be a regular $16$-gon inscribed into ${\bf E}^2\cap S^4$
and let ${\PP}_3$ be a set of $16$ pairwise antipodal points on
${\bf E}^3\cap S^4$ with covering radius $r_c=33.547...^{\circ} $. For the details of the construction 
of ${\PP}_3$ see \cite{gft}. Finally, let
${\PP}$ be the convex hull of ${\PP}_2\cup {\PP}_3$. If $F$ is any
facet of ${\PP},$ then it is easy to see that $F$ is a 4-dimensional
simplex having two vertices in ${\PP}_2$ and three vertices in ${\PP}_3$.
If one projects $F$ from the center $o$ of $S^4$ onto $S^4,$ then the
projection $F'$ is a $4$-dimensional spherical symplex. Among its five vertices
there are two vertices say, $a$ and $b$ lying in ${\bf P}_2$. Here $a$ and $b$ must be consecutive
vertices of the regular $16$-gon inscribed into ${\bf E}^2\cap S^4$, while the remaining
three vertices must form a triangle inscribed into ${\bf E}^3\cap S^4$ with circumscribed circle ${\cal C}$ of radius $r_c$. Now, Let $c'$ be the center of ${\cal C}$ and $c$ be an arbitrary point of ${\cal C}$ moreover, let
$m$ be the midpoint of spherical segment $ab$. Clearly,  $am=11.25^{\circ }$, $c'm=90^{\circ }$ and $cc'=r_c$ on $S^4$. If $s$ denotes the center of the circumscribed sphere of $F'$ in $S^4$, then $s$ is a point of the spherical segment $c'm$. Let $as=bs=cs=x$, $sm=y$ and $c's=90^{\circ }-y$.
Now, the cosine theorem applied to the spherical right triangles $\Delta ams$ and $\Delta cc's$ implies that
$$\cos x=\cos 11.25^{\circ } \cdot \cos y, \text{\ {\rm and}\ } \cos x=\cos r_c \sin y.$$
By solving these equations for $x$ and $y$ we get that
$x=50.572...^{\circ }<$ $\arccos \sqrt{\frac{2}{5}}$ $=50.768...^{\circ }$ finishing the proof of Theorem~\ref{negyedik} for $d=5$.

For $d=6$ we need to construct $64$ pairwise antipodal points on $S^5$ with
covering radius at most $\arccos \sqrt{\frac{5}{12}}=49.797...^{\circ }$. 
In order to achieve this let us take two $3$-dimensional subspaces
${\bf E}^3_1$ and ${\bf E}^3_2$ in ${\bf E}^6$ such that they are totally 
orthogonal to each other. For $i=1,2$ let ${\PP}_i$ 
be a set of $32$ pairwise antipodal points on
${\bf E}^3_i\cap S^5$ with covering radius $r_c=22.690...^{\circ} $. 
For the details of the construction 
of ${\PP}_i, i=1,2$ see \cite{gft}. Finally, let
${\PP}$ be the convex hull of ${\PP}_1\cup {\PP}_2$. If $F$ is any
facet of ${\PP}$, then it is easy to see that $F$ is a 5-dimensional
simplex having three vertices both in ${\PP}_1$ and in ${\PP}_2$.
If one projects $F$ from the center $o$ of $S^5$ onto $S^5,$ then the
projection $F'$ is a 5-dimensional spherical symplex. It follows from the construction above that two
spherical triangles formed by the two proper triplets of the vertices of $F'$ have circumscribed circles
${\cal C}_1$ and ${\cal C}_1$ of radius $r_c$. If $c_i$ denotes the center of ${\cal C}_i, i=1,2$ and 
$s$ denotes the center of the circumscribed sphere of $F'$ in $S^5$, then it is easy to show that $s$ is in fact, the midpoint of the spherical segment $c_1c_2$ whose spherical length is of $90^{\circ }$. Thus, if $x$ denotes the spherical radius of the circumscribed sphere of $F'$ in $S^5$, then the cosine theorem applied to the proper spherical right triangle implies that
$$\cos x=\cos r_c \cdot \cos 45 ^{\circ }  .$$
Hence, it follows that
$x=49.278...^{\circ }< \arccos \sqrt{\frac{5}{12}}=49.797...^{\circ }$ finishing the proof of Theorem~\ref{negyedik} for $d=6$.

\kkk

Finally, we state the following immediate corollary, which proves the Illumination Conjecture for convex bodies of constant width of dimensions $4,5$ and $6$.

\begin{cor}
If $\WW$ is a convex body of constant width in ${\bf E}^4$, then $I(\WW)\leq 12$. Moreover, if $\WW$ is a convex body of constant width in ${\bf E}^d$ with $d=5,6$, then $X({\WW})\leq 2^d$. 
\end{cor}

\section{Weakly neighbourly antipodal con\-vex po\-ly\-topes}
\label{three}

As we shall see in this section the X-ray Conjecture is naturally connected to a special kind of antipodality that raises a challanging new question for further study. In 
order to phrase it properly we need to introduce the notion of a {\it weakly neighbourly antipodal convex polytope} along
with a simple statement whose straightforward proof we leave to the reader. (Actually, neither antipodality nor weak neighbourleness are new notions however, it seems that this is the first time when they are considered simultaneously.)

\begin{definition}
The convex polytope $\PP \subset {\bf E}^d$ is called a weakly neighbourly convex polytope if any two vertices of $\PP$ lie on a face of $\PP$. Moreover, the convex polytope $\PP \subset {\bf E}^d$ is called an antipodal convex polytope if any two vertices of $\PP$ lie on parallel supporting hyperplanes.
\end{definition}

\begin{proposition} 
Let $\PP \subset {\bf E}^d$ be a $d$-dimensional weakly neighbourly antipodal convex polytope. If the number
of the vertices of $\PP$ is $v$, then $v\le X(\PP)$.
\end{proposition}

Now, it is clear that the X-ray Conjecture implicitly contains the following also quite challanging conjecture
that is worth phrasing independently.

\begin{con}
If $\PP \subset {\bf E}^d$ is an arbitrary $d$-dimensional weakly neighbourly antipodal convex polytope, then the number of vertices of $\PP$ is at most $3\cdot 2^{d-2}$.
\end{con}

On the one hand, it has been known for quite some time  (see \cite{dg}) that if $\PP \subset {\bf E}^d$ is an antipodal convex polytope,
then the number of vertices of $\PP$ is at most $2^d$. On the other hand, based on the very recent paper \cite{ss} one can easily verify that the above conjecture holds in dimensions $2$ and $3$.

\vspace{1cm}

\medskip

\noindent
K\'aroly Bezdek,
Department of Mathematics and Statistics,
2500 University drive N.W.,
University of Calgary, AB, Canada, T2N 1N4.
\newline
{\sf e-mail: bezdek@math.ucalgary.ca}
 
\smallskip

\noindent
Gy\"orgy Kiss,
Department of Geometry, Mathematical Institute,
E\"otv\"os Lo\-r\'and University,
P\'azm\'any P\'eter s\'et\'any 1/c,
H-1117 Budapest, Hungary.
\newline
{\sf e-mail: kissgy@cs.elte.hu}

\end{document}